\documentclass[a4paper,leqno,11pt,twoside]{amsart}
\usepackage{amsmath,amsfonts,amssymb,amsthm,amscd}
\usepackage[utf8]{inputenc}
\usepackage[english]{babel}
\usepackage[colorlinks=true,citecolor=blue, urlcolor=blue, linkcolor=blue,pagebackref]{hyperref}
\usepackage[top=1.2in,bottom=1.2in,left=1.in,right=1.in]{geometry}

\usepackage[all, cmtip]{xy}
\usepackage{enumerate}
\usepackage{cancel}

\usepackage{tikz}
\usepackage{tikz-cd}
\usepackage{graphicx}
\usepackage{pgf,tikz}
\usetikzlibrary{arrows}

\newtheorem{thm}{Theorem}[section]
\newtheorem{lem}[thm]{Lemma}
\newtheorem{cor}[thm]{Corollary}
\newtheorem{prop}[thm]{Proposition}

\theoremstyle{definition}
\newtheorem{rem}[thm]{Remark}

\newcommand{\ZZ}{{\mathbb Z}}

\newcommand{\cO}{{\mathcal O}}

\newcommand{\cA}{{\mathcal A}}

\newcommand{\cP}{{\mathcal P}}

\newcommand{\cRH}{{\mathcal {RH}}}

\newcommand{\tC}{{\widetilde C}}

\newcommand{\lra}{\longrightarrow}
\newcommand{\ra}{\rightarrow}

\newcommand{\PP}{{\mathbb P}}

\DeclareMathOperator{\Ker}{Ker}
\DeclareMathOperator{\Ext}{Ext}
\DeclareMathOperator{\End}{End}
\DeclareMathOperator{\im}{Im}

\DeclareMathOperator{\Div}{Div}
\DeclareMathOperator{\Hom}{Hom}
\DeclareMathOperator{\Pic}{Pic}

\newcommand{\restr}[1]          {\vert_{#1}}
\title[Simplicity of some Jacobians]{Simplicity of some Jacobians with many automorphisms}

\author{Juan Carlos Naranjo$^{1,2}$}
\address{Juan Carlos Naranjo \newline 1. Departament de Matem\`atiques i Inform\`atica, Universitat de Barcelona, Gran Via de les Corts Catalanes, 585, 08007 Barcelona, Spain \newline 2. Centre de Recerca Matemàtica, Edifici C, Campus Bellaterra, 08193 Bellaterra, Spain }
 \email{jcnaranjo@ub.edu}

\author{Angela Ortega}
\address{Angela Ortega \newline 
 Institut f\"ur Mathematik, Humboldt Universit\"at zu Berlin, Unter den Linden 6, 10099 Berlin,  Germany}
\email{ortega@math.hu-berlin.de}

 \author{Gian Pietro Pirola}
 \address{Gian Pietro Pirola, Universit\`a degli Studi di Pavia, Dipartimento di Matematica, Via Ferrata 5,  27100 Pavia, Italy  }
 \email{gianpietro.pirola@unipv.it}

\author{Irene Spelta}
\address{Irene Spelta \newline 
 Centre de Recerca Matemàtica, Edifici C, Campus Bellaterra, 08193 Bellaterra, Spain }
\email{ispelta@crm.cat}

\begin{document}
\begin{abstract}
We study an explicit $(2g-1)$-dimensional family of Jacobian varieties of dimension $\frac{d-1}2(g-1)$, arising from quotient curves of unramified cyclic coverings of prime degree $d$ of hyperelliptic curves of genus $g\ge 2$.  By using a deformation argument, we prove that the generic element of the family is simple. Furthermore, we completely describe their endomorphism algebra, and we show that they admit a rank $\frac{d-1}2-1$ group of non-polarized automorphisms. As an application of these results, we
prove the generic injectivity of the Prym map  for \'etale cyclic coverings of hyperelliptic curves of odd prime degree under some slight numerical restrictions. This result generalizes in several directions previous results on genus 2. 
\end{abstract}

\maketitle
\section{Introduction}

Explicit high-dimensional families of Jacobians of curves with many non-polarized automorphisms are difficult to find. Many natural examples come from curves which map to elliptic curves. Hence, the corresponding Jacobians are strongly non-simple. 
The other few explicit examples we are aware of are discussed in \cite{ell}: they are obtained from branched covers of the projective line, and they admit real multiplication by a subfield of even index at most 10 in a primitive cyclotomic field. In \cite{ell} (see reference therein) it is explained how, via this construction, one can obtain almost all the other examples  discovered by Brumer, Mestre, and Tautz, Top, and Verberkmoes. The geometry of such examples is so rich that they have been exploited in other interesting contexts (\cite{albano pirola, lo16, i}). 

In this paper, we deal with a higher dimensional generalization of one of the examples presented in \cite{ell}. 
We discuss the geometric properties of a  family of Jacobian varieties having a collection of automorphisms not preserving the polarization. Our aim is to prove that the generic element of the family is simple by using deformation and monodromy arguments. Moreover, we fully describe the endomorphism algebras of these Jacobians.
Our argument fits in the research line, dating back to Ekedhal and Serre, which investigates the decomposability of abelian varieties. In the last years, several developments have been achieved, mainly based on group algebra or/and period matrix techniques (see for instance \cite{rr} and references therein). In this sense, our result provides a completely new criterion for the simplicity of certain Jacobians. 

\bigskip

 More precisely, we consider a cyclic unramified covering $f:\tilde C\rightarrow C$ of an irreducible complex smooth hyperelliptic genus $g\ge 2$ curve, of degree an odd prime number $d=2k+1$. Let $\sigma$ be the automorphism on $\tC$ such that $C= \tC /\langle \sigma \rangle$.
The hyperelliptic involution  $\iota$ on $C$ lifts to an involution $j$ on $\tilde C$. We studied the family of the quotient curves 
(and their Jacobians) $C_0:=\tilde C/\langle j \rangle$ of genus $g(C_0)=k(g-1)$ depending on $2g-1$ parameters, the dimension of the moduli space $\mathcal H_g$ of hyperelliptic curves of genus $g$. 

 Since the quotient map $\pi_0: \tilde C\ra C_0$ is ramified, the pullback map $\pi_0^*$ is injective. Thus, we identify $JC_0$ with its image in $J\tilde C$. It is known (\cite{ries}) that the automorphisms $\sigma^i+\sigma^{-i}$, for $i=1,\ldots, k$, of $J\tilde C$ restrict to automorphisms of  $JC_0$, denoted by $\beta_i$. It turns out that these automorphisms of $JC_0$ do not preserve the polarization.   

 Our main Theorem is the following:

 \begin{thm}\label{teo simplicity}
 For a generic  cyclic unramified covering $f:\tilde C\rightarrow C$ of degree $d$, the Jacobian $JC_0$ is simple and the endomorphism algebra $End_{\mathbb Q}(JC_0)$ of $JC_0$ is the totally real field   $\mathbb Q(\xi+\xi^{-1})$, where $\xi$ is a primitive $d$-th root of the unity.  Moreover, the group of automorphisms has rank $k-1$ and is generated by any $k-1$ elements in $\{\beta_1,\ldots ,\beta_{k}\}$.
 \end{thm}
The proof of this result is based on a deformation argument which rely on Hodge theory and a monodromy argument. Indeed, we show that the infinitesimal deformations of 
$(\tilde C, C)$ which preserve the hyperelliptic structure, ``rigidify'' the cotangent space to $JC_0$. In this way, we deduce the simplicity, as well as the description of the endomorphism algebra. Thus, by applying the Dirichlet's unity theorem, we obtain the description of the automorphism group.

Theorem \ref{teo simplicity} has a direct application to Prym theory, allowing to generalize the main result of \cite{nos} in several directions.  In order to state our result, we need to introduce some notation. Let $\cRH_{g}[d] $ denote the coarse moduli space of cyclic unramified degree $d$ covers of hyperelliptic genus $g$ curves. As before, $d$ is an odd prime number.  
It is well-known that every element $f:\tilde C\ra C$ in $\cRH_{g}[d] $ is described by a pair $(C, \langle \eta \rangle )$, where $C$ is a hyperelliptic genus $g$ curve and $\eta\in JC\smallsetminus \{0\}$ is a non-trivial $d$-torsion point. The Prym variety $P(\tilde C, C)$ of the covering is defined as the connected component containing the origin of the kernel of the Norm map $J\tilde C\ra JC$ and thus, $\dim P(\tilde C, C)=(d-1)(g-1)$. Moreover, the principal polarization on $J\tilde C$ induces a polarization $\tau$ on it of type $\delta = (1, \dots, 1, d, \dots,d)$, with $d$ repeated $(g-1)$-times. We study the generic injectivity of the Prym map:
 \[
 \cP_{g}[d]: \cRH_{g}[d]\ra \cA_{(g-1)(d-1)}^{\delta }
 \]
sending the isomorphism class of the covering  $f:\tilde C\ra C$ to the isomorphism class of $(P(\tilde C, C), \tau)$.

The Prym map $\cP_{g}[d]$ has already been considered in the literature. The case $d=2$ has been treated by Mumford in \cite{mu}, where it is shown that $\cRH_{g}[2]$ consists of several irreducible components.
In \cite{albano pirola}, the authors showed that the generic fiber is positive dimensional for $d=3$ and $g\leq 4$, as well as $d=5$ and $g=2$. When $d=7$ and $g=2$ the map is dominant onto its image, and it is generically finite of degree 10, see \cite{lo16}. Quite recently, the infinitesimal Torelli theorem stated for $d\geq 6$ and $g=2$ in \cite{ag}, 
has been improved in \cite{nos} by proving the generic Torelli theorem for $d$ prime and such that $\frac{d-1}{2}$ is prime too (i.e. the so-called Sophie Germain primes). 
 Roughly speaking, the proof in \cite{nos} is based on the existence of a (non-polarized) isomorphism $P\ra JC_0\times JC_0$, where $C_0$ is, as before, the quotient of $\tilde C$ by a lifting $j$ of the hyperelliptic involution on $C$. By using the endomorphisms on $P$ and the polarization $\tau $, the curve $C_0$ and the automorphisms $\beta_i$ of the Jacobian $JC_0$, are  intrinsically  recovered from $(P,\tau)$.  This part uses
 the simplicity of $JC_0$ and some properties of the endomorphism ring $\End(JC_0)$, which are proved under the hypotheses $g=2$ and $\frac {d-1}2$ prime.  
 The last part of the proof relies on a ``theta-duality" procedure inspired by the properties of the Fourier-Mukai transform. Using Theorem \ref{teo simplicity}, and  modifying  the strategy in \cite{nos},  we obtain the following  generalization of
 the main Theorem in \cite{nos}:

\begin{thm} \label{teo Prym}
Assume that $d$ is a prime number with $ (d-1)(g-1)\ge 7$, $d\ge 3$, and  that $g\ge 2$ is not congruent with $3$ modulo $d$, then the Prym map $\cP_{g}[d]$ is generically injective. 
\end{thm}

This improves the Theorem in \cite{nos}  in two directions: it removes the assumption of $\frac {d-1}2$ being prime in genus $2$, and it extends to hyperelliptic curves of genus $g\ge 3$ (under some slight numerical constrains).
One can also remove the congruence condition if $g>2d-1$ (see Remark \ref{high_genus}). 
Furthermore, under more general assumption, we prove the following:

\begin{thm} \label{gen_finiteness}
For $g\ge 3,\ d\ge 5$, and for $g=2, \ d\ge 6$ the map is generically finite onto its image.
\end{thm}

The generic finiteness of the Prym map was proven for $g=2$ in \cite{ag}. The argument for $g\ge 3$ in Section 3 is independent of the rest of the paper. 

\vskip 5mm
\textbf{Acknowledgements:} The first and the fourth authors were partially supported by  the Departament de Recerca i Universitats de la Generalitat de Catalunya (2021 SGR 00697), the Proyecto de Investigaci\'on PID2023-147642NB-100 and the Spanish State Research Agency, through the Severo Ochoa and María de Maeztu Program for Centers and Units of Excellence in R\&D (CEX2020-001084-M). The third author is partially supported by PRIN project {\em Moduli spaces and special varieties} (2022). The third and the fourth authors are  members of GNSAGA (INdAM). The fourth author was partially supported by INdAM-GNSAGA project CUP E55F22000270001.

  \section{The simplicty of $JC_0$}
In this section, we introduce some notation and we prove Theorem \ref{teo simplicity}. In the rest of the paper $d$ is an odd prime number and we set $k:=\frac{d-1}{2} $. Let $f: \widetilde C\ra C$ be a cyclic \'etale cover of a hyperelliptic curve $C$ of genus $g$. Let $\eta\in JC$ be a non-trivial $d$-torsion point such that $\langle \eta \rangle $ equals the kernel of $f^*: JC\ra J\widetilde C$. 
Let $\mathcal{RH}_g[d]$ be the coarse moduli space of cyclic \'etale covers of hyperelliptic curves of degree $d$. In order to use a monodromy argument, we will need the irreducibility of $\mathcal{RH}_g[d]$ and of the moduli space $\widetilde{\mathcal{RH}}_g[d]$ of classes of pairs $[C,\rho]$, where $\rho \in JC[d]\smallsetminus \{0\}$. Notice that this is not true for $d=2$ by \cite[section 7]{mu} (see also \cite{BF} and \cite{naranjo}).

\begin{prop}\label{irr}
The moduli spaces $\mathcal{RH}_g[d]$ and $\widetilde{\mathcal{RH}}_g[d]$ are irreducible for $d\ge 3$.
\end{prop}
\begin{proof}
The irreducibility of the moduli space $\mathcal{RH}_g[d]$ can be deduced from the proof of Theorem 1 in \cite{BF}, which is proven for $m$-gonal curves with $m\geq 3$ and cyclic étale coverings of degree $n \geq 2$, but it can be easily extended to hyperelliptic curves ($m=2$) and $n$ odd. Following the notation of \cite{BF}, it is enough to show that the conjugacy class of $(\alpha_{ij}; (i \ j))$ contains the element $({\bf 0}; (1 \ 2) ) \in 
\ZZ_n^2 \rtimes S_2$. Setting $\tau=((-\alpha, 0, \ldots, \alpha, 0, \ldots, 0); (1 \ 2) )$, one checks that $\tau \cdot (\alpha_{ij}; (i \ j)) \cdot \tau^{-1} = ({\bf 0}; (1 \ 2) )$.  For the sake of clarity, we include a more intuitive argument by induction of the irreducibility. Observe that there is a natural finite unramified map $\widetilde{\mathcal{RH}}_g[d] \rightarrow \mathcal{RH}_g[d],$ hence it is enough to prove the irreducibility of $\widetilde{\mathcal{RH}}_g[d]$. On the other hand, there is nothing to prove for $g=1,2$ since in this case all curves are hyperelliptic. 
We work in the Deligne-Mumford compactification  $\overline{\mathcal {H}}_{g,d} $ of $\widetilde {\mathcal{RH}}_g[d]$ which has a finite map on the closure of the hyperelliptic locus $\overline {\mathcal H}_g$ in $\overline {\mathcal M}_g$.  Let us consider a hyperelliptic curve $H_0$ of genus $g-2$ with two marked Weierstrass points $w,w' \in H_0$. Let $E_1, E_2$ be two elliptic curves with marked origins $O, O'$. By gluing $w$ with $O$ and $w'$ with $O'$ we obtain a semistable hyperelliptic curve $H\in \overline{\mathcal H}_g$  whose normalization is the disjoint union $H_0\sqcup E_1 \sqcup E_2$. Observe that $JH=JH_0 \times E_1 \times E_2$. Let $(a,b,c) \in JH[d]$ be any $d$-torsion point in $JH$. We want to connect $(H,(a,b,c))\in \overline{\mathcal {H}}_{g,d} $ with $(H,(0,p,0))$, where $p$ is a fixed  $d$-torsion point in $E_1$. First, we observe that fixing $H_0, E_2$ we can use the irreducibility in genus $1$ to connect $(0,p,0)$ with $(0,b,0)$. If $a=c=0$, we are done, otherwise we  use the irreducibility of $\overline{\mathcal{H}}_{g-1,d}$ (fixing $E_1$) to connect $(a,b,c)$ with $(a',b,0)$. Now fixing $E_2$ and applying again induction we can connect $(a',b,0)$ with $(0,b,0)$ and we are done. 
\end{proof}

Let $\sigma$ be a generator of the Galois group of $f$ and $\iota$ the hyperelliptic involution of $C$. We have $\langle \sigma \rangle\cong \mathbb{Z}/d$. We denote $\sigma$ as well the induced automorphism on $J\widetilde C$. 
Since  the covering $\widetilde C$ is defined as $\text{Spec}(\mathcal{O}_C\oplus \eta\oplus \dots\oplus \eta^{d-1})$ and $\iota^*\eta^i= \eta^{d-i}$, the hyperelliptic involution of $C$ lifts to an involution $j$ on $\widetilde C$. The resulting covering $\widetilde C\ra \mathbb{P}^1$ is Galois with Galois group the dihedral group $D_d=\langle j,\sigma \mid j^2=\sigma^d=1, \; j\sigma j=\sigma^{-1}\rangle$. Notice that all the involutions $j\sigma^i\in D_d$, with $i=0, \dots, d-1$, yield intermediate quotient curves $C_i:= \widetilde C/\langle j\sigma ^i\rangle $. Since all the  involutions are  conjugated, the corresponding curves  $C_i$ are all isomorphic. Hence, from now on we focus on $C_0$.

Let $\chi$ be a generating character of $\langle\sigma\rangle$, i.e. $\chi(\sigma)=\xi$, where $\xi$ is a primitive $d$-th root of the unity that we fix from now on. By means of the projection formula, we obtain the decomposition
\begin{equation}\label{suma_diretta}
H^0(\widetilde C, \omega_{\widetilde C})=\bigoplus_{i=0}^{d-1}H^0(C, \omega_{C}\otimes \eta^i),
\end{equation}
where $H^0(C, \omega_{C}\otimes \eta^i)$ is the eigenspace associated to the character $\chi^i$.  By the Riemann-Roch formula, we have that $h^0(C, \omega_{C}\otimes \eta^i)=g-1$, for $i>0$. Notice that the $g$-dimensional vector space $H^0(C, \omega_{C})$ corresponds to the eigenspace with eigenvalue 1, thus it is the only one which is invariant under the action of $\sigma$ on $H^0(\widetilde C, \omega_{\widetilde C})$. 

Let $P=P(\widetilde C, C)$ be the Prym variety of the covering $f$. It is an abelian variety of dimension $(d-1)(g-1)$ and polarization of type $(1, \ldots, 1,d, \ldots,d)$ with $d$ appearing $g-1$ times. We have that the dual of the tangent space of $P$ at $0$ decomposes as follows:
\begin{equation}\label{tang Prym}
	T_0P^*=\bigoplus_{i=1}^{d-1}H^0(C, \omega_{C}\otimes \eta^i).
\end{equation}

Now we describe the tangent space to $JC_0$. Let $s$ be a section in $H^0(C, \omega_{C}\otimes \eta^i)\subset H^0(\tC, \omega_{\tC})$. We have that $\sigma(s)=\xi^is$. Furthermore, the relations in the dihedral group imply the equalities \begin{equation}\label{group law dihedral}
	\sigma j (s)=j\sigma^{-1}(s)=j(\xi^{-i}s)= \xi^{-i}j(s).
\end{equation}
Hence, if $s\in H^0(C, \omega_{C}\otimes \eta^i)$, then $js\in H^0(C, \omega_{C}\otimes \eta^{-i}) $. Let us define 
the vector spaces 
\[
V(i):=\{ s+js \mid s\in H^0(C, \omega_{C}\otimes \eta^i)\} \subset H^0(C, \omega_{C}\otimes \eta^i)\oplus H^0(C, \omega_{C}\otimes \eta^{-i}),
\] 
for $i=1,\ldots,k$.

We have the following: 
\begin{prop}\label{cotang JC_0}
	The cotangent space of $JC_0$ at the origin is 
	\begin{equation}\label{tang JC0}
		T_0JC_0^*=\bigoplus_{i=1}^k V(i).
	\end{equation}
\begin{proof}
This follows from $H^0(C_0, \omega_{C_0} )=H^0(\tilde C, \omega_{\tilde C})^{\langle j \rangle}$, and the equalities (\ref{suma_diretta}) and \eqref{group law dihedral}. 
\end{proof}
\end{prop}

\begin{rem} \label{action_beta}
	Let us observe that \[(\sigma^i+\sigma^{-i})(s+js)=(1+j)(\sigma^i+\sigma^{-i})(s) \quad \forall \; s\in H^0(C, \omega_{C}\otimes \eta^l). \]
	It follows that the action of $\beta_i=(\sigma^i+\sigma^{-i})_{\vert JC_0}$ on $T_0JC_0^*$ is diagonalizable with eigenspaces $V(l)$ and corresponding eigenvalues $\xi^{il}+\xi^{-il}$, for $l=1,\dots,k$.
\end{rem}
 Now, let $\epsilon$ be in $H^1(C, T_C)^{\langle\iota\rangle}$, that is, an infinitesimal deformation of $C$ preserving the hyperelliptic structure. Deformations as such come from deformations of $(C, \langle\eta\rangle)$, in other words, we see these deformations as deformations of $\tilde C$. Therefore, we have that $H^1(C, T_C)^{\langle\iota\rangle}\subset H^1(\tilde C, T_{\tilde C})$. The equivariance of the multiplication map yields the following: 
 \begin{prop}
     The multiplication maps 
     \[
     \begin{aligned}
     m_i: H^1(C, T_C)\cdot (H^0(C, \omega_C \otimes \eta^i)\oplus H^0(C, \omega_C \otimes \eta^{-i})) & \ra  H^1(C, \eta^i)\oplus H^1(C, \eta^{-i}) \\ &\cong 
     (H^0(C, \omega_C \otimes \eta^{-i})\oplus H^0(C, \omega_C \otimes \eta^{i}))^*
     \end{aligned}
     \]
      restrict as follows: \begin{equation}\label{multiplication V(i)}
         H^1(C, T_C)^{\langle\iota\rangle}\cdot V(i)\ra V(i)^*.
     \end{equation}
 \end{prop}

In particular, we consider the ``diagonal" multiplication maps $\overline m_i:V(i) \rightarrow H^0(C, \omega_C^2)^{\langle\iota\rangle}$ sending $s+js$ to $s\cdot js$, $i=1, \dots, k$. Let $L_i$ be the image of $\overline m_i$. We have the following:
  
\begin{lem}\label{lem torsion}
	 Let  $C$ be a general hyperelliptic curve. Then $ L_i \cap  L_l= \{0\}$ for all $i \neq l$.
	\begin{proof}
		By contradiction, let us assume that $L_i \cap L_l\neq \{0\}$ for some $i\neq l$. 
Then we have $$s\cdot js= t\cdot jt$$
for some non-trivial $s\in H^0(C, \omega_C\otimes\eta^i)$ and $t\in H^0(C, \omega_C\otimes\eta^l)$. We set $s'=js$ and $t'=jt$ and we denote by 
$D(s)= p_1 +\ldots+p_{2g-2}$ and $D(t)=q_1 +\ldots +q_{2g-2}$  the divisor of $s$ and $t$ respectively. Under our assumption,  we obtain an equality of divisors:\[
p_1 +\ldots+ p_{2g-2} + \iota(p_1) + \ldots+ \iota(p_{2g-2}) = q_1 + \ldots+ q_{2g-2} + \iota(q_1) +\ldots+ \iota(q_{2g-2}).\]
Now we set \begin{equation}\label{decompositon div}
    D(s) =A +\iota B \quad\text{and}\quad D(s')=B + \iota A, \; \text{with}\; A,B \in \Div(C)
\end{equation} 
in such a way that
$A + B = D(t)$. Changing $\eta$ with $\eta^{-1}$ if necessary, we may assume that $\deg B\geq g-1$. Thus, \[\mathcal{O}_C(D(t)-D(s'))=\mathcal{O}_C(A-\iota A)= \eta^{l-i}\neq \mathcal O_C,\]
with $\deg A \leq g - 1$. 
Now, we consider a map $\gamma $ defined on the symmetric product $C^{(g-1)}$ as follows:
   \begin{align*}
   \gamma: C^{(g-1)}\ra JC\\
    A\mapsto A- \iota A.
\end{align*}
Since $\gamma(A)=A- \iota A= 2A-K_C$, we obtain that the divisor $\Gamma:=\gamma(C^{(g-1)})$ belongs to $\vert 2\Theta\vert$. Thus, it suffices to show that $2\Theta\cap JC[d]$ is empty for a general $C$. We claim that for any symmetric theta divisor $\Theta$, $2\Theta\cap JC[d] \neq \emptyset$
implies that $\Theta\cap JC[2d]\neq \emptyset$. Indeed, if
$D-\iota D \in 2\Theta\cap JC[d]$, then $dD \sim d\iota D$. For
a theta characteristic $\kappa$ corresponding to the symmetric 
theta divisor $\Theta$ we have the equality in $JC$
$$
d(D-\kappa) = d \iota D - d \iota\kappa = \iota(dD-d\kappa) = -d(D-\kappa).
$$
Hence $2d(D-\kappa)=0$ and $D-\kappa \in \Theta \cap JC[2d]$, which proves the claim.

Now, let us assume that there exists a $2d$-torsion point in $\Theta\cap JC[2d]$. Since the moduli space $\widetilde{\mathcal{RH}}_g[d]$ is irreducible (see Proposition \ref{irr}),  and since $C$ is general, using monodromy, we obtain that the whole group $JC[2d]$ of $2d$-torsion points is contained in $\Gamma$. Now, by \cite[Theorem C]{par}, there exists a bound on the number of $2d$-torsion point on the theta divisor of a principally polarized abelian variety. Since this bound is strictly smaller than the cardinality of $JC[2d]$, we obtain a contradiction. 
 \end{proof}
\end{lem} 
In order to better describe the multiplication map in  \eqref{multiplication V(i)}, we first analyze how we can exhibit special deformations lying in $H^1(C, T_C)^{\langle \iota\rangle}$. We have the following:
\begin{lem}\label{Schiffer}
    Let $D\in \Div(C)$ be an effective divisor (of degree less or equal than $2g-2$) and assume that $D=F+A+\iota A+W$, where $W$ is a combination of Weierstrass points and $F$ does not contain divisors of the form $x+\iota(x)$, for $x\in C$. Then the intersection  
    \[H^0(C, T_C(D)\restr{D})\cap H^1(C,T_C)^{\langle \iota \rangle }\subset H^1(C,T_C)\]
    is generated by the Schiffer deformations $\{\epsilon_{a_i}+\epsilon_{\iota(a_i)}, \epsilon_{w_l}\}_{i,l}$, where $A=\sum a_i$ and $W=\sum w_l$. 
    In particular, $\dim H^0(C, T_C(D)\restr{D})\cap H^1(C,T_C)^{\langle \iota \rangle }\leq \deg (A)+\deg (W)$.
    \begin{proof}
        From the exact sequence \begin{equation}\label{seq tangs}
        0\ra T_C\ra T_C(D)\ra T_C(D)\restr{D}\ra 0,         
        \end{equation}
        we have that $H^0(C, T_C(D)\restr{D})$ injects inside $H^1(C, T_C)$ and, by definition, the image is generated by the Schiffer deformation $\epsilon_p$, for $p\in D$. By considering the ones which are left invariant by the action of $\iota$, we get that the intersection $H^0(C, T_C(D)\restr{D})\cap H^1(C,T_C)^{\langle \iota \rangle }$ is as claimed. The statement on the dimension follows easily by recalling that we are admitting divisors $D$ with non-reduced support. 
    \end{proof}
\end{lem}

This allows us to state the following:
\begin{lem}\label{lem surj}
    For every $s\in H^0(C, \omega_C\otimes \eta^i)$ the map $H^1(C,T_C)^{\langle \iota \rangle }\ra V(i)^*$ sending $\epsilon\mapsto \epsilon s+ \epsilon js$ is surjective. 
    \begin{proof}
        Let $D$ be the divisor which is the intersection of the divisors attached to $s$ and $js$. As in the previous Lemma, we put $D=F+A+\iota A+W$ and we consider the following exact sequence:
        \[0\ra \mathcal{O}_C(D)\xrightarrow{(\cdot s, \cdot js)} (\omega_C\otimes \eta^i)\oplus ( \omega_C\otimes \eta^{-i})\ra \omega_C^2(-D)\ra 0. \]
        After tensoring with $T_C$ and combining with the exact sequence in \eqref{seq tangs}, in cohomology we have:
        \[
        \xymatrix{  &  & 0\ar[d] &  & \\
 &  & H^0(C, T_C(D)\restr{D})\ar[d] & & \\
 &  & H^1(C, T_C) \ar[d]\ar[dr]^\alpha & & \\
0\ar[r] & H^0(C, \omega_C(-D))\ar[r] & H^1(C, T_C(D)) \ar[r] \ar[d] & H^1(C, \eta^i)\oplus H^1(C, \eta^{-i})\ar[r] & ... \\
&  & 0 &  &. }
\]
We claim that $\dim (\ker \alpha) \cap H^1(C, T_C)^{\langle\iota\rangle}\leq g$. Indeed, 
$(\ker \alpha) \cap H^1(C, T_C)^{\langle\iota\rangle}$ applies to $H^0(C,\omega_C(-D))$ and therefore there is an exact sequence:
\[
0 \ra H^0(C, T_C(D)\restr{D})^{\langle\iota\rangle} \ra \ker \alpha\cap H^1(C, T_C)^{\langle\iota\rangle}
\ra H^0(C,\omega_C(-D)) \ra \ldots
\]
By Lemma \ref{Schiffer}, we know that $\dim H^0(C, T_C(D)\restr{D})^{\langle\iota\rangle}\leq \deg (A)+\deg (W)$. On the other hand, since $h^0(C, \omega_C(-D))\leq g-(\deg (A)+\deg (W))$, we have: 
\[
\dim \ker \alpha\cap H^1(C, T_C)^{\langle\iota\rangle}\leq \dim H^0(C, T_C(D)\restr{D})^{\langle\iota\rangle} + H^0(C,\omega_C(-D))\leq g.
\]
Since the image of $\alpha $ restricted to the invariant part $H^1(C, T_C)^{\langle\iota\rangle}$ (of dimension $2g-1$) is necessarily a subset of $V(i)^*$, which has dimension $g-1$,
 we have shown the surjectivity.   
       
    \end{proof}
\end{lem}

From now on, we assume that  $(C, \langle\eta\rangle)$ is general, and we focus on the simplicity of $JC_0$. 
Let us assume that there exists an abelian subvariety $A\subset JC_0$, we shall show that it is trivial, that is, either $A= \{0\}$ or $A= JC_0$. 
Let  $W:=T_0 A^*$.
By using the polarizations on $JC_0$ and $A$, we can think of  $W$ as a subspace of $\oplus _{i=1}^k V(i)$. Let $\epsilon$ be in $H^1(C, T_C)^{\langle\iota\rangle}$.
By \eqref{multiplication V(i)}, we have $\epsilon \cdot \bigoplus V(i)\subset \bigoplus V(i)^*  $. In the same way, since we assume that $A$ is an abelian subvariety of $JC_0$ deforming with $(C, \langle\eta\rangle)$, we have that $\epsilon\cdot W\subset W^*.$ 
Using the existence of a polarization, we identify $W^*$ with $\overline W$. First, we observe the following: 
\begin{prop}\label{part 2 intersection}
    If $W\cap V(i)\neq \{0\}$, then $V(i)\subseteq W$.
    \begin{proof}
        Indeed, let $s_i+js_i$ be an element in this intersection. By Lemma \ref{lem surj}, the map \[H^1(C,T_C)^{\langle \iota \rangle }\ra V(i)^*\] sending $\epsilon\mapsto \epsilon s_i+ \epsilon js_i$ is surjective. 
    Therefore,  for every $ t_i+jt_i\in V(i)^*$, there exists $\epsilon\in H^1(C, T_C)^{  \langle\iota\rangle}$ such that  $\epsilon (s_i+js_i)=t_i+jt_i$. By the invariance of $W$ by deformations and conjugation we obtain  $V(i)\subseteq W$, as desired.
    \end{proof}
\end{prop}

Now, let $0\neq \Omega=\sum \Omega_i$ be an element of $W$. By the equality (\ref{tang JC0}), we can decompose $\Omega=\sum_{i=1}^k \Omega_i $, with $\Omega_i=s_i+js_i\in V(i), s_i\in H^0(C, \omega_{C}\otimes \eta^i)$. We have the following ``killing Lemma''.

\begin{lem}\label{lemma killer}
   Let $1\leq i,l\leq k$ be to different indices, and let $\Omega_i\in V(i), \Omega_l \in V(l)$ be two non-trivial elements. Then, there exists $\epsilon \in H^1(C, T_C)^{\langle\iota\rangle}$ such that $\epsilon \Omega_i=0$, and  $\epsilon(\Omega_l)\neq 0$.
    \begin{proof}
  Let $s_i \in H^0(C, \omega_{C}\otimes \eta^i)$ be such that $\Omega _i=s_i+js_i$. The vanishing of $\epsilon \Omega_i$ for some $\epsilon \in H^1(C, T_C)^{\langle\iota\rangle}$ is equivalent to the conditions $\epsilon s_i=\epsilon js_i=0$. Indeed, $\epsilon s_i$ belongs to $H^1(C, \eta^i)$ while $\epsilon js_i$ to $H^1(C, \eta^{-i})$. Since $\epsilon $ and $j$ commute, it is enogh to look for a $\epsilon$ such that $\epsilon s_i=0$.

Furthermore, let us observe that we can apply a deformation $\epsilon \in H^1(C,T_C)^{\langle \iota \rangle}$ to the initial $\Omega_i, \Omega_l$ and then conjugate to obtain new elements in the same subspaces. By Lemma \ref{lem surj}, applied to the section $s_i$, we can then assume that $s_i$ is ``generic'', namely we have that $R=(s_i)_0$ does not contain Weierstrass points, neither divisors of the form $x+\iota x$, for $x\in C$. Then, we consider the short exact sequence attached to $s_i$:
  	\[
    0\ra \mathcal{O}_C \stackrel{\cdot s_i}{\ra}\omega_{C}\otimes \eta^i\ra \omega_{C}\otimes \eta^i\restr{R}\ra 0.
    \]
     Tensoring with $T_C$, we obtain 
  \[
  0\ra T_C \ra \eta^i\ra \eta^i_{\restr{R}}\ra 0,
  \]
 that, in cohomology, yields 
  \begin{equation}\label{mult}
      0\ra \mathbb{C}^{2g-2}\ra H^1(C, T_C)\xrightarrow{\alpha}H^1(C, \eta^i)\ra 0.
  \end{equation}
Notice that $\alpha (\epsilon)=\epsilon s_i \in H^0(C, \omega_C \otimes \eta ^{-i})^*\cong H^1(C, \eta^i)$.  By Grassmann's formula, we have:
\begin{equation} \label{dim_ker_alpha}
\dim \Ker(\alpha)\cap H^1(C, T_C)^{\langle\iota\rangle} \ge 2g-2+ \dim H^1(C, T_C)^{\langle\iota\rangle} -\dim H^1(C,T_C)=g.
\end{equation}
  
We put now  $\Omega_l=s_l+js_l \neq 0$. Let $D$ be the divisor which is the intersection of the divisors attached to $s_i$ and $s_l$. Recall that, under our assumption of generality, it does not contain Weierstrass points, neither divisors of the form $x+\iota x$. We consider the sequence:
\[0\ra T_C(D) \ra \eta^i\oplus \eta^{l}\ra \omega_C\otimes \eta ^{i+l}(-D)\ra 0, \]
where the first map is given by the multiplications with the sections $s_i,s_l$ on each summand and the second sends $(x,y)$ to $s_lx-s_iy$.
We can assume $\deg D\leq  g-2$. Indeed, if $\deg D\geq g$ we can do the same argument replacing $s_i$ with $js_i$, while the case $\deg D= g-1$ is excluded by Lemma \ref{lem torsion}, since in this case, up to a constant, $s_i j(s_i)=s_lj(s_l)$.  Hence, we consider the following diagram in cohomology: 
 \[
        \xymatrix{  &  & 0\ar[d] &  & \\
 &  & H^0(C, T_C(D)\restr{D})\ar[d] & & \\
 &  & H^1(C, T_C) \ar[d]\ar[dr]^\beta & & \\
0\ar[r] & H^0(C, \omega_C\otimes \eta^{i+l}(-D))\ar[r] & H^1(C, T_C(D)) \ar[r] \ar[d] & H^1(C, \eta^i)\oplus H^1(C, \eta^{l})\ar[r] & ... \\
&  & 0 &  &. }
\]
In order to prove our statement, we need to show that there exists $\epsilon\in H^1(C, T_C)^{\langle\iota\rangle}$ such that $\epsilon\cdot s_i=0$ and $\epsilon \cdot s_l\neq 0$ for $l\neq i$ and $\Omega_l\neq 0$. In other words, we search for $\epsilon\in \ker(\alpha)\cap  H^1(C, T_C)^{\langle\iota\rangle}$ such that $\beta(\epsilon)\neq 0$. As in Lemma \ref{lem surj}, we have 
$$\dim \ker(\beta)\cap H^1(C, T_C)^{\langle\iota\rangle}\leq \dim H^0(C, T_C(D)\restr{D})^{\langle\iota\rangle} + H^0(C,\omega_C\otimes \eta^{i+l}(-D)). $$
Under our assumption of generality on $D$, by Lemma \ref{Schiffer}, we have $\dim H^0(C, T_C(D)\restr{D})^{\langle\iota\rangle}=0$, while $h^0(C,\omega_C\otimes \eta^{i+l}(-D))\leq g-1 $. Using \eqref{dim_ker_alpha}, we conclude.

 \end{proof}
\end{lem}
Now we are ready to describe completely the cotangent space $W$. 
\begin{prop}\label{cotang W}
	The following equality holds: $W=\bigoplus_{i\in I} V(i)$, for some $I\subset \{1,\dots, k\}$.
	\begin{proof}

  Let $\Omega=\sum_{i\in I} \Omega_i \in W$ for a certain subset of indices $I\subset \{1,\dots, k\}$, and assume that $\Omega_i\neq 0$ for $i\in I$. We can assume without lost of generality that $I=\{1,\ldots,l\}$, $1\le l \le k$. 

  We claim that $W\cap V(i)\neq \{0\}$ for all $i=1,\ldots, l$. If $l=1$ there is nothing to prove, hence we assume $l\ge 2$. We will prove that $W\cap V(1)\neq \{0\}$, the same argument works for the other indices.  By Lemma \ref{lemma killer}, there exists a deformation $\epsilon \in H^1(C, T_C)^{\langle\iota\rangle}$ sending to zero $\Omega_2$ and such that $\epsilon \Omega_1$ is not zero.  This means that $\epsilon\Omega= \sum_{i\in I\setminus \{2\}}\epsilon \Omega_i$. As already remarked, we have $\epsilon\Omega\in \overline{W}$, thus $\overline{\epsilon\Omega}\in {W}.$
 
 Notice that 
 $\epsilon\Omega_i=\epsilon s_i+\epsilon js_i$ belongs to $\overline{V(i)}$. Hence, using conjugation again, we obtain that $\overline{\epsilon\Omega}_i$ belongs to $V(i)$.  
  Now we start again with the new element $\overline{\epsilon\Omega}$ which has, at least, one index less. Iterating we obtain an element in $W$ which belongs to $V(1)$, this proves the claim. 
  
  All in all, using Proposition \ref{part 2 intersection}, we obtain that: 
  \[
  V(1)\oplus \ldots \oplus V(l) \subseteq W. 
  \]
  If the inclusion were not an equality, we could choose  a new $\Omega$ in $W\setminus \oplus_{i=1}^lV(i)$ and start again.
  Since $\dim W$ is finite, this procedure ends and we get the statement.  
\end{proof}
\end{prop}

Now we are ready for the proof of the Theorem \ref{teo simplicity}. Let $\widetilde C\ra C$ be a general element in $\mathcal{RH}_g[d]$.  We have the following:
\begin{thm}\label{simplicity}
	The Jacobian $JC_0$ is simple.
	\begin{proof}
    We proceed by contradiction. Let us assume that there exists a non-trivial abelian subvariety $A\subset JC_0$ and set $B=JC_0/A$. By Proposition \ref{cotang W}, we have that $W:=T_0A^*=\oplus_{i\in I}V(i)$, for a set of indices $I$ that we can assume to be $I= \{1,\dots, l\}$.  
	
	Observe that the case $l=1$ is not possible. Indeed, this would give $T_0A^*=V(1)$. By Remark \ref{action_beta}, the automorphism $\beta_1$ leaves invariant $A$ and the restriction $\beta_{1 \restr{A}}$ would determine an endomorphism of $A$ which is equal to $(\xi+\xi^{-1})\mathrm{Id}$. Since $(\xi+\xi^{-1})\in \mathbb R$ is not an integer, this is impossible. 
		Hence, let us assume $l\geq 2$. We work now in the irreducible moduli space $\widetilde {\mathcal{RH}}_g[d]$ of classes of pairs $[C,\rho]$, where $\rho \in JC[d]\smallsetminus \{0\}$ (see Proposition \ref{irr}).

  The isogeny decomposition of $JC_0$ must be preserved under monodromy, that is, it remains the same when we replace $[C, \eta]\in \widetilde {\mathcal{RH}}_g[d] $ with $[C, \eta^l]$ for any $l=1,\dots, d-1$.  
		
	 The action $\eta\mapsto \eta^l$ is given by an automorphism belonging to the Galois group $G:=Gal(\mathbb{Q}(\xi)/\mathbb{Q})$. Since $d$ is prime, we have $G\cong \mathbb{Z}/(d-1)\mathbb{Z}$.

	 Since, by definition $V(i)=V(d-i)$, we take into account only the automorphisms  
     $\eta\mapsto \eta^l$ for $l\in \{2,\ldots , k\}$. 
	
	 In order to investigate the possible configurations for $W$, we have to study the Galois orbits of the units $\xi, \dots, \xi^{k-1}, \xi^k$. Indeed, these orbits precisely provide the only ones possible partitions of $\{1, \dots, k\}$ that remain fixed under the action of the automorphisms of $G$. 
	 More precisely, let $\eta\mapsto\eta^l$ be one of these automorphisms. If $l$ had order equal to $k$ or $2k$, then the automorphism would act as a permutation of the set 
     $\{1,\dots, k\}$. This would say that there is a single Galois orbit, namely 
     $W=\oplus_{i=1}^k V(i)$. In other words, that $A=JC_0$, and this is incompatible with 
    our assumptions. Therefore, let us assume that $m:=ord(l)<k$. Let $I_1\cup \dots \cup 
 I_{\frac{k}{m}}$ be a partition of the set $\{1\, \dots, k\}$ given by the Galois orbits 
     of $\xi, \dots, \xi^k$ under the automorphism  $\eta\mapsto\eta^l$. Notice that all
     the Galois orbits are of length equal to $m$.  Then $W=\oplus_{i\in I_r} V(i)$ for a certain $I_r$ among $I_1, \dots, I_{\frac{k}{m}}$. We consider the endomorphism $\sum_{i\in I_r} \beta_i$. By construction, it restricts to an endomorphism of $W$ which equals $\lambda\text{Id}$, with $\lambda=(\sum_{i\in I_r}\xi^i+\xi^{-i})$. Again, since $\lambda \in \mathbb R$ is not an integer, we obtain a contradiction. 
		 
	\end{proof}
\end{thm}
 Let $\text{End}^0(JC_0):=\text{End}(JC_0)\otimes \mathbb{Q} $ be the endomorphism algebra of $JC_0$. For simplicity, we set $D:=\text{End}^0(JC_0)$ and we let $Z$ be the centre. Furthermore, being $d\mapsto d'$ the positive Rosati involution on $D$, we denote by $Z_0$ the fixed field of such involution restricted to $Z$. Since $JC_0$ is simple, by Albert's classification, $D$ is a division algebra with $Z_0\subseteq Z\subseteq D$ satisfying one of the following:
 \begin{itemize}
     \item[1.] $Z_0=Z=D$ is a totally real field;
     \item[2.] $D$ is a quaternion algebra (totally definite/totally indefinite) over a totally real field $Z=Z_0$;
     \item[3.] $D$ is a division algebra over the CM-field $Z$ (over the totally real subfield $Z_0$). 
 \end{itemize}

Under our assumptions of genericity, the deformation argument exploited in Proposition \ref{cotang W} allows us to exclude the last two possibilities. Indeed, we have the following:
 \begin{thm}\label{end}
     The endomorphism algebra $D$ of $JC_0$ is the totally real field   $D=\mathbb Q(\xi+\xi^{-1})$.
     \begin{proof}
      We claim that it is enough to show that $\mathbb Q(\xi+\xi^{-1})$ is the maximal field contained in $D$. Indeed, this would automatically exclude cases 2 and 3, and the maximality would yield $D$ as in the statement. 

      In order to prove the claim, let us assume by contradiction that there exists an intermediate field extension $\mathbb Q(\xi+\xi^{-1})\subsetneq L\subseteq D$. By the primitive element theorem, there exists a $\alpha$ such that $L=\mathbb Q(\alpha)$.  Since $\mathbb Q(\xi+\xi^{-1})\neq L$, we have that \begin{equation}\label{degree}
          [L:\mathbb Q]>[\mathbb Q(\xi+\xi^{-1}):\mathbb Q]=k.
      \end{equation}  By an abuse of notation, we will denote with the same $\alpha$ the endomorphism induced on $T_0JC_0^*$. This is diagonalizable: indeed, its minimal polynomial divides the polynomial of $\alpha$ over $\mathbb Q$, which is irreducible and separable. Therefore, let $\bigoplus_{l=1}^rW_l$ be the eigenspace decomposition of $\alpha$. We apply the same deformation argument used in Proposition \ref{cotang W}: we know that each of the subspaces of $H^{1,0}(JC_0)$ that remains invariant under deformations by $\epsilon\in H^1(C, T_C)$ is either one of the $V(i)$, $i=1,\ldots, k$, or a direct sum of several of them. Thus, we get $W_l= \oplus_{i\in I_l} V(i)$ for some $I_l\subset \{1, \dots, k\}. $  This yields a contradiction: indeed, by \eqref{degree}, we get $r>k$, which is impossible since there are exactly $k$ subspaces $V(i)$. 
     \end{proof}
 \end{thm}

\begin{cor}\label{aut}
    The group of automorphisms of $JC_0$ has rank $k-1$ and is generated by any subset of cardinality $k-1$ of $\{\beta_1,\ldots ,\beta_k\}$.

\begin{proof}
    The statement about the rank is a consequence of the Dirichlet Theorem on units of a field. The second claim follows from $\beta_1+\ldots +\beta_k=1$, which is straightforward by replacing $\beta_i$ with $\sigma ^i+ \sigma ^{-i}$.
\end{proof}

\end{cor}

Observe that Theorems \ref{simplicity}, \ref{end}, together with Corollary \ref{aut}, give Theorem \ref{teo simplicity} in the Introduction.

\section{Generic finiteness of the Prym map}

The aim of this section is to prove that the Prym map 
 \[
 \cP_{g}[d]:\cRH_{g}[d]\lra \cA_{(g-1)(d-1)}^{\delta},
 \]
is generic finite as soon as the dimensions of the source and the target allow it. This is equivalent to the 
surjectivity of the codifferential map 
 \[
 d\cP_{g}^*[d]: T^*_{[P,\tau]}\cA_{(g-1)(d-1)}^{\delta} \lra T^*_{[C, \langle \eta \rangle]} \cRH_{g}[d] ,
 \]
at the general points $[C, \langle \eta \rangle]\in \cRH_{g}[d]$ and $[P,\tau]$ in the image of the Prym map. Notice that, combining \cite[Remark 2.5]{ag} together with the isomorphism $T^*_{[C, \langle \eta \rangle]} \cRH_{g}[d]\cong H^0(C, \omega_C^2)^{\langle \iota \rangle}$, we can equivalently describe $d\cP_{g}^*[d]$ as a multiplication map followed by a projection as follows
\begin{equation}\label{codiff prym hyp}
 \bigoplus_{i=1}^k H^0(C, \omega_C\otimes \eta^i) \otimes H^0(C, \omega_C\otimes \eta^{-i})\stackrel{m}{\ra}  H^0(C, \omega_C^2)  \stackrel{p}{\ra}  H^0(C, \omega_C^2)^{\langle \iota \rangle}.
\end{equation}

We start by showing the following:
\begin{prop}\label{surjectivity_multiplication} Let $L$ be the line bundle corresponding to the hyperelliptic linear series on $C$. We assume that $g\geq 5$. Let $\eta \in Pic^0(C)$ such that  $ \eta^2$ is not of the form $L^2(-D)$, where $D$ is an effective divisor of degree $4$. Then,  the map $$m_{\eta }:H^0(C,\omega_C\otimes \eta )\otimes H^0(C,\omega_C\otimes \eta ^{-1})\to H^0(C,\omega_C^2)$$
is surjective.\end{prop}

\begin{rem}\label{remark} The assumption on $\eta $ in Proposition \ref{surjectivity_multiplication} includes both,
$\eta \neq 0$ and that $\omega_C\otimes \eta$  
is  generated by its sections. 
Indeed, this last statement is equivalent to say that $\eta\neq \cO(p-q)$ for any two points $p,q$, since if there exists a base point $p$ for $\omega_C\otimes \eta$, then $h^0(C,\omega_C\otimes \eta (-p))=g-1$, that is $h^0(\eta^{-1}(p))=1.$ 
Let us see that our hypothesis imply that $\eta$ is not of the form  $\cO(p-q)$: if so,
 since the hyperelliptic involution $\iota$ acts as $-\mathrm{Id}$  on $\Pic^0(C)$, we would also have $\eta = \cO(\iota q-\iota p)$. This implies $\eta^2=\cO_C(p+\iota q-q-\iota p)=\cO_C( p+\iota p+q+\iota q-2q-2\iota p))=L^2(-2q-2\iota p)$, which contradicts the assumption.
\end{rem}

\begin{proof}
Assume by contradiction that $m_{\eta }$ is not surjective.
The dual of the multiplication map is given by: 
\[
\varphi: H^1(C,T_C) \ra \Hom^s(H^0(C,\omega_C \otimes \eta), \ H^0(C,\omega_C \otimes \eta^{-1})^*),
\]
which can be described as follows: let $0\neq \epsilon \in H^1(C,T_C)=\Ext^1(\omega_C,\mathcal O_C)$ corresponding to a non-split extension:
\begin{equation}\label{basic_extension}
0\ra \mathcal O_C \ra E_{\epsilon} \ra \omega_C \ra 0. 
\end{equation}
We tensor this sequence with $\eta $:
\[
0\ra \eta  \ra E_{\epsilon}\otimes \eta \ra \omega_C \otimes \eta\ra 0,
\]
and we consider the coboundary map:
\[
\partial_{\epsilon} : H^0(C,\omega_C\otimes \eta ) \ra H^1(C,\eta) \cong H^0(C,\omega_C\otimes \eta^{-1})^*
\]
Then, for any $s\in H^0(C,\omega_C\otimes \eta)$, $\varphi(\epsilon )(s)=\partial_{\epsilon} (s)$.
Therefore, the non-injectivity of the multiplication map implies the existence of a non-split extension $[E_{\epsilon}]$ as above, 
such that $\partial_{\epsilon} (s)=0$ for all $s$. 
This is equivalent to the condition:  $h^0(C,E_{\epsilon} \otimes \eta )=g-1$. We use the following result by Segre-Nagata and Ghione (see \cite[Example 7.2.14]{lazars}):

\begin{lem}\label{lemma Ghio}
	Let $C$ be a curve and $E$ a rank $v$, degree $p$ vector bundle on it. Then there exists a subline bundle $M\subset E$ such that \[\deg(M)\geq \frac{p}{v}-\frac{v-1}{v}(g-1).\]
\end{lem}
In our case we have $v=2$ and $p=2g-2$, thus the Lemma yields the existence of a subline bundle $M$ such that $\deg (M)\geq \frac{g-1}{2}$. Since by assumptions, $g\geq 5$, we have $\deg M\geq 2.$ So, up to saturation (which could increase the degree of $M$), 
we can assume the existence of a sequence of locally free sheaves:
\begin{equation}
0\to M\to E_{\epsilon}\otimes \eta \to  N \to 0 \label{ass},\end{equation}
and a diagram:
\begin{equation}
 \label{diagramma croce}
		\begin{tikzcd}
&	&0	\arrow[d] & &\\
&	&M \arrow{dr}{\tau}\arrow[d] & &\\
0\arrow[r]&\eta  \arrow[r]& E_{\epsilon}\otimes \eta 
\arrow[r]\arrow[d]&\omega_C\otimes\eta \arrow[r]&0\\
& &N	\arrow[d] & &\\
& &0 & & 
		\end{tikzcd}
		\end{equation}
Notice that	the map $\tau$ is necessarily non-zero, since otherwise there would be an injective map $M\ra \eta $ which is impossible since the degree of $M$ is positive.
Observe that $N\cong \omega_C \otimes \eta ^2 \otimes M^{-1}$ and that the vertical exact sequence implies that:
\[
h^0(N)+h^0(M)\geq g-1, \qquad \deg N +\deg M=2g-2.
\]
Now we discuss the possible values for the degrees and the dimensions of the space of sections of $N$ and $ M$.  We prove first several claims:

\textbf{Claim 1:} $\deg(M)<2g-2$, and therefore $\deg(N)>0$.

Indeed, otherwise $\tau $ is a non-trivial map between two line bundles of the same degree, hence $M=\omega _C\otimes \eta $, $N=\eta$ and the map $\tau$ must be the multiplication by a non-zero constant. The inverse of $\tau $ would give a section of $E_{\epsilon}\otimes \eta \ra \omega_C \otimes \eta$ and therefore the extension (\ref{basic_extension}) would split. 

\textbf{Claim 2}\label{claim 2} : $1\leq h^0(C,N)\leq g-2 $, which yields $h^0(C, M)\ge 1.$

Indeed, if $h^0(N)= g-1$, by Clifford inequality, 
$g-1 \leq \frac{\deg N}{2} + 1 \leq g-2  +1$,  then $N$ is a power of the line bundle $L$ corresponding to the hyperelliptic linear series: $N= L^{(g-2) }$. Thus $ N= \omega_C \otimes L^{-1}$ and $M= L\otimes\eta^2$. In particular $\deg(M)=2$. By Lemma \ref{lemma Ghio}, $\deg M\geq \frac{g-1}{2}$, this is only possible for $g=5$. In this case, we can assume $h^0(C, L\otimes\eta^2)=0$, otherwise we would have $\eta^2=\mathcal{O}_C(p-q)$, for certain $p,q\in C$, which is against our assumption.  Thus $E_{\epsilon}\otimes \eta$ fits in the short exact sequence
$$
0\to L \otimes\eta^2 \to E_{\epsilon}\otimes \eta \to  L^{3} \to 0 $$
 with zero coboundary map $H^0(L^3) \to H^1(L\otimes\eta^2)$.  
Equivalently, the multiplication map 
$$
 H^0(L^3) \otimes H^0(L^3\otimes \eta^{-2}) \to H^0(L^6\otimes \eta^{-2})
$$
 is not surjective. Note that $L^3\otimes\eta^{-2}$ is base point free by our hypothesis. Since both sides of the map are of the same dimension, the base-point-free pencil trick shows that the multiplication map is an isomorphism, which gives 
 a contradiction.

Moreover, $h^0(C,N)$ must be positive since otherwise $h^0(C,M)=g-1$ and then the map $\tau $ induces an isomorphism 
$H^0(C,M)\ra H^0(C,\omega_C\otimes \eta )$. Since $\tau $ can be seen as a section of $\omega_C\otimes \eta \otimes M^{-1}$, this section can not have zeros, then $M\cong \omega_C \otimes \eta$ and, as before,  
the initial extension  (\ref{basic_extension}) splits.

\textbf{Claim 3:} $\deg N>1$.

Otherwise, by claim 2 we would have that $h^0(C,N)=1$, so $N$ is either trivial or  $N=\cO_C(p)$ for some point $p\in C$. The trivial case  implies $M\cong \omega_C \otimes \eta^2$ and, as in the proof of Claim 1, this yields a contradiction to the non-splitness of the extension  (\ref{basic_extension}).
Assume that $N=\cO_C(p)$, then  $M=\omega_C\otimes \eta ^2(-p)$. The inclusion $\omega_C\otimes \eta ^2(-p)\hookrightarrow \omega_C\otimes \eta$ implies $ \eta^2(-p)\hookrightarrow \eta$ and thus $ \eta^{-1}(p)$ has a section, so $\eta$ is of the form $\cO_C(p-q)$, a contradiction.

Summarizing, we have found the following restrictions on the degrees on $M,N$ (with $M\otimes N=\omega_C \otimes \eta^2$):
\[
\deg(N)+\deg(M)=2g-2, \ \ \deg(M)\geq 2, \ \ \deg(N)\geq 2,
\]
and on the dimensions of the spaces of global sections:
\[
h^0(C,N)+h^0(C,M)\geq g-1,\ \  1\leq h^0(C,N)\leq g-2, \ \ 1\leq h^0(C,M)\leq g-1. 
\]

\begin{lem}
With these restrictions, we have one of the following possibilities:
\begin{enumerate}
\item $M$ and $N$ are special, that is $h^1(C,M)>0$ and $h^1(C,N)>0,$ 

\item $N=L$ and $M=\omega_C \otimes \eta^2 \otimes L^{-1}$,
\item $M=L$ and $N=\omega_C \otimes \eta^2 \otimes L^{-1}$.
\end{enumerate}
\end{lem}
\begin{proof}    
If one of them, let's say $N$, is not special, then 
$h^0(C,\omega_C\otimes N^{-1})=0 $ and therefore $N$, which can be represented by an effective divisor,
 must have degree $g+k$ with $k\geq 0$. Notice that   $g+k\leq 2g-4$,  that is $k\leq g-4$. By the Riemann-Roch Theorem, we have that $h^0(C,N)=k+1$.
 On the other hand, $M$ has degree $g-2-k$ and  then it is special. By the Clifford inequality $h^0(C,M)\leq \frac{\deg M}2 +1$. Then
$$g-1\leq h^0(N)+h^0(M)\leq  
k+1+\frac{\deg M}2 +1=k+2+\frac{g-2-k}2.$$
This implies that $ k\geq g-4$ and $M$ satisfies the equality in Clifford's Theorem, hence $M=L$ and $N=\omega_C \otimes \eta ^2 \otimes L ^{-1}$. The same proof applies verbatim if $M$ is not special.
\end{proof}

Now we can finish the proof of Proposition \ref{surjectivity_multiplication}. Assuming that the multiplication map is not surjective we have found the existence of the line bundles $M,N$ 
satisfying a series of numerical constrains. Our aim now is to prove that these conditions lead to a contradiction with the hypothesis of the Proposition. To do so, we analyze the three possibilities that appear in the previous Lemma. 

Let us  first assume that we are in condition $(i)$, namely that $M, N$ are special. Then, since $C$ is hyperelliptic, we get that $M=L^r(D_1)$ and $N=L^s(D_2)$ where $L$ is, as before, the $g_2^1$ on $C$, and 
$D_1$ and $D_2$ are effective. Moreover, we have $h^0(M)=r+1$ and
$h^0(N)=s+1$ with $r+s\geq g-3.$ Then, we have:
\[
\omega_C \otimes \eta ^2=M\otimes N = L^{r+s}(D_1+D_2)=L^{g-3}(E),
\]
where $E$ is an effective divisor of degree $4$. Since $\omega_C=L^{g-1}$, we obtain that
$L^2\otimes \eta^2=\cO_C(E)=L^4(-\iota E) $. 
Therefore, $\eta ^2= L^2(-\iota E)$ which contradicts our hypothesis.

Easily $(ii)$ yields a contradiction. Indeed,   if $N=L$ and $M=\omega_C\otimes \eta^2\otimes L^{-1}$, the map $\tau $ gives a section of $L\otimes \eta^{-1}$. Then $\eta $ is of the form $L(-p-q)=\cO_C(\iota(p)-q)$. 

 Thus, it only remains to handle condition $(iii)$, i.e., $M=L$. Let us first observe that, as in Claim 2, we have to study the combination $M=L, N=\omega_C\otimes \eta^2\otimes L^{-1}=L^3\otimes \eta^2$ for $g=5$. The strategy is the same as in Claim  2. Indeed, we can assume 
that $L^3\otimes \eta^2$ is not special, in other words that $E_{\epsilon}\otimes \eta$ fits in the short exact sequence
$$
0\to L \to E_{\epsilon}\otimes \eta \to  L^3\otimes \eta^2 \to 0 $$
with zero coboundary map $H^0(L^3\otimes \eta^2 ) \to H^1(L)$.  


Equivalently, the multiplication map 
$$
 H^0(L^3) \otimes H^0(L^3\otimes \eta^{2}) \to H^0(L^6\otimes \eta^{2})
$$
 is not surjective. Note that $L^3\otimes\eta^{2}$ is a base point free pencil by our hypothesis. Since both sides of the map are of the same dimension, the base-point-free pencil trick shows that the multiplication map is an isomorphism, which gives 
 a contradiction.
\end{proof}

\begin{rem}\label{torsion-points}
Observe that the general element $(C, \eta) \in \cRH_g[d]$ satisfies the condition of the Proposition \ref{surjectivity_multiplication}. 
Suppose $g\geq 5$ is odd, then $\kappa=L^{\frac{g-1}{2}}$ is a theta characteristic and let 
$\Theta_{\kappa}$ the corresponding symmetric theta divisor. One verifies that  $C^{(2)} - C^{(2)} \subset \Theta_{\kappa}$.  
So, if for all $\eta \in JC[d]$, $\eta = \cO(2L-D) \in C^{(2)} - C^{(2)}$, (with $D\in C^{(4)}$), then all the $d$-torsion points would be contained in $\Theta_{\kappa}$, which contradicts the bound in \cite{par}. In case  $g>5$ is even, one can use the same argument with the theta characteristic 
$\kappa=L^{\frac{g-2}{2}}\otimes \cO(w)$,
where $w$ is a Weierstrass point.
\end{rem}

Proposition \ref{surjectivity_multiplication} is no longer true for genus $g\leq 3$ by dimension reasons. Neither in $g=4$, where, interesting enough, the image of the multiplication map is of codimension one. So, in order to get a surjective map onto $H^0(C, \omega_C)^{\langle \iota\rangle}$ one needs to involve another summand of  \eqref{codiff prym hyp}. For this reason, we will treat cases $g=3,4$ with direct methods.

\begin{prop} \label{mult_small_genus}
Let $g=3$ or 4 and $d\geq 5$. For any $d$-torsion element $ \eta \in JC[d]\smallsetminus \{ 0\}$ the map $p\circ m$
is surjective.
\end{prop}
\begin{proof}
Since $d\geq 5$, it is enough to check that the restriction of $p\circ m$ to the first two summands: 
\begin{equation}\label{mult map g4}
    H^0(C,\omega_C\otimes \eta )\otimes H^0(C,\omega_C\otimes \eta ^{-1})\oplus H^0(C,\omega_C\otimes \eta^2 )\otimes H^0(C,\omega_C\otimes \eta ^{-2}) \to H^0(\omega_C^2)^{\langle \iota \rangle}.
\end{equation}
is surjective.

{\bf Case $g=4, d\geq 5$}. Consider the surjective difference map:  
\[ 
C^{(2)}\times C^{(2)} \to JC, \quad (z,w)\mapsto [z-w]. 
\] 
So $\eta=D_1-D_2$, for some divisors $D_i\in C^{(2)}$. Therefore, setting $D_2=x+y, $ for $x, y\in C$, we get $\omega_C\otimes \eta= L^3(D_1-D_2)= L(D_1+\iota x+\iota y)$. Thus, we have the inclusion $ H^0(C, L)\subset H^0(C,\omega_C\otimes \eta )$ which yields 
an isomorphism
\[
H^0(C,\omega_C\otimes \eta )\cong H^0(C, L)\oplus \langle u\rangle,\]
where $\langle u\rangle=H^0(C, \cO_C(\iota x+\iota y))$. Similarly, we have \[H^0(C,\omega_C\otimes \eta ^{-1})\cong H^0(C, L)\oplus \langle ju\rangle.\] 
Hence, under the above isomorphisms, the multiplication map 
\[m_{\eta }:H^0(C,\omega_C\otimes \eta )\otimes H^0(C,\omega_C\otimes \eta ^{-1})\to H^0(C,\omega_C^2 )\] 
can be identified with the multiplication map:
\[(H^0(C, L)\oplus \langle u\rangle)\otimes (H^0(C, L)\oplus \langle ju\rangle)\to H^0(C,\omega_C^2 )\] 
which clearly contains $\Lambda^2 H^0(C,L)$
in its kernel.  Since the source and the target of $m_{\eta }$ have the same dimension,
$m_{\eta }$ is not surjective. 
We set $H^0(C, L)=\langle s_1, s_2\rangle$. It is easy to show that $u\cdot ju\in \langle s_1^2, s_2^2, s_1s_2\rangle$. 
Thus, taking the second summand in \eqref{mult map g4} and considering the analogous decompositions $$H^0(C, \omega \otimes \eta^2)\cong H^0(C,L)\oplus \langle u '\rangle, \quad  H^0(C, \omega \otimes \eta^{-2})\cong H^0(C,L)\oplus \langle ju '\rangle, $$ we prove that the following elements
\[ 
s_1^2, \ s_2^2, \ s_1s_2, \ s_1u+s_1ju, \ s_2u+s_2ju, \ s_1u'+s_1ju', \ s_2u'+s_2ju' \]
are linear independent in $H^0(C, \omega_C^2)^{\langle \iota \rangle}$ and can be used to generate $\im p\circ m$.  Therefore, since $h^0(C, \omega_C^2)^{\langle \iota\rangle}=7$, we can conclude the surjectivity of the map. The check on the dimension is done using computer calculations in SAGE (https://amor.cms.hu-berlin.de/~speltair/SageCasesg3g4.html).

{\bf Case $g=3, d\geq 5$}. 
 Notice that in this case $h^0(C,\omega_C\otimes\eta^i)=2$ for every $i$, while $h^0(C, \omega_C^2)^{\langle \iota\rangle}=5$. 
 Moreover, by the base-point-free-pencil trick, the  restricted maps $m_{\eta}$ and $m_{\eta^2}$ are injective. We fix the following bases:\[H^0(C,\omega_C\otimes \eta )=\langle t_1,t_2\rangle, \quad H^0(C,\omega_C\otimes \eta^2 )=\langle w_1,w_2\rangle.\]
Again by means of computer calculations in SAGE (see the link above), we can show that the following elements in $H^0(C, \omega_C^2)^{\langle \iota \rangle}$
\[
t_1jt_1, \ t_2jt_2, \ w_1jw_1,\  w_2jw_2, \ t_1jt_2+t_2jt_1
\] 
are linear independent. Since by definition they belong to $\operatorname{Im}d\mathcal{P}^*$, the map $p\circ m$ is surjective. 
\end{proof}
\begin{thm}
     Assume that either $g=2$ and $d\geq 6$ or $g\geq 3$ and $d\geq 5$, or 
     $g\ge 5$ and $d\geq 3$. Then the Prym map $\cP_{g}[d]$ is generically finite.
\end{thm}
 \begin{proof}
   For $g=2$ this is the main theorem in the paper \cite{ag} by Agostini. 
      
        The generic finiteness of $\cP_{g}[d]$  is equivalent to  the surjectivity of the composition map \eqref{codiff prym hyp} for a general element $[C, \langle \eta \rangle] \in \cRH_{g}[d]$. 
        Since the projetion is surjective  it is enough to show that the 
        multiplication map $m$ is surjective.
Proposition  \ref{surjectivity_multiplication}
together with Remark \ref{torsion-points} show that $m$ is surjective for $g\geq 5$. The cases 
$g=3,4$ are solved in Proposition \ref{mult_small_genus}.

\end{proof}



\section{Application to the generic injectivity of the Prym map}

In this section we will use Theorem \ref{simplicity} to show that, under some numerical conditions, the Prym map 
 \[
 \cP_{g}[d]:\cRH_{g}[d]\lra \cA_{(g-1)(d-1)}^{\delta},
 \]
 is generically injective. First, we collect some properties of $JC_0$.
 \begin{prop}[Ries \cite{ries}, Ortega \cite{ortega}]\label{Properties}
 	 The following properties hold.
 	\begin{itemize}
 		\item $P=\text{Im}(1-\sigma)$ and $\sigma^*\tau=\tau$. Thus $\sigma$ restricts to an automorphism of $P$ as polarized abelian variety.
 		\item For $i\in \{1,\dots, d-1\}$, we have $JC_i=\im (1+j\sigma^i) \subset P$. Moreover, the equality $\sigma(1+j\sigma^i)=(1+j\sigma^{i-2})\sigma$, gives  that $\sigma$ maps $ JC_i$ to $ JC_{i-2}$. In particular, $\sigma^i: JC_0\xrightarrow{\cong} JC_{d-2i}$.
 		\item For $i \in \{1,\dots, k\}$, the maps $\sigma^i+\sigma^{-i}$ restrict to automorphisms of $JC_0$ (not preserving its principal polarization), we call $\beta_i$ to the restriction. Moreover, $\beta_i'=\beta_i$, i.e., it is symmetric with respect to the Rosati involution.
 		\item For $i \in \{1,\dots, k\}$, the map $\psi_i: JC_0\times JC_0\ra P$,  $(x,y)\mapsto x+\sigma^i(y)$, is an isomorphism such that \begin{enumerate}
 			\item $\psi_i^*(\tau)=\begin{pmatrix}
 				2\lambda_0 & \lambda_0\beta_i\\
 				\lambda_0\beta_i &2\lambda_0
 			\end{pmatrix}$, where $\lambda_0: JC_0\ra \widehat{JC}_0$ is the natural principal polarization;
 		\item $\sigma^i\circ \psi_i= \psi_i\circ \begin{pmatrix}
 			0 &-1\\
 			1&\beta_i
 		\end{pmatrix}$.
 		\end{enumerate}
 	\end{itemize}
 \end{prop}

As in \cite{nos}, we factor the Prym map as the composition of two maps: \begin{equation}
		 \cP_{g}[d]: \cRH_{g}[d]\xrightarrow{\cP_1[d]} \mathcal{D}_k \xrightarrow{\cP_2[d]}\cA_{(g-1)(d-1)}^{\delta },
\end{equation}
where $\mathcal{D}_k$ parametrizes classes of objects $(C_0, \beta_1, \dots, \beta_k)$, $\cP_1[d] $ sends $(C, \eta)$ to $(C_0, \beta_1, \dots, \beta_k)$ and $\cP_2[d] $ sends $(C_0, \beta_1, \dots, \beta_k)$ to $(JC_0\times JC_0, \begin{pmatrix} 
		2\lambda_0 & \lambda_0\beta_i\\
		\lambda_0\beta_i &2\lambda_0
\end{pmatrix}) $ which, by Proposition \ref{Properties}, lies in $\text{Im} \cP_{g}[d]$. 

First, we prove that $\cP_2[d]$ is generically injective. Let $(P, \tau)$ be a general element in $\text{Im} \cP_{g}[d]$. Then, by \cite[Proposition 3.1]{nos}, we have that \begin{equation}\label{autom P}
	\{\alpha\in \text{Aut}(P,\tau): \alpha(x)=x, \, \forall x\in K(\tau)\}\cong \langle \sigma \rangle.
\end{equation}
This is shown by using the assumptions on $\alpha$ to construct an $\tilde \alpha$ acting on $J\widetilde C$, then proving that $\tilde \alpha$ is compatible with the principal polarization on $J\tilde C$ and that it acts as the identity on $f^*JC$. 

Hence, for the general $(P, \tau)\in \text{Im} \cP_{g}[d] $,  $\langle\sigma\rangle$ is uniquely determined. Therefore, we now work with the triplet $(P, \tau,\langle \sigma \rangle)$. 

The following Proposition is proved in \cite{nos} by means of arithmetic arguments that require the extra conditions $g=2$ and $k$ prime. We will use instead Theorem \ref{teo simplicity} to prove this more general version. The strategy is the same, we recall the main facts by the convenience of the reader. 

\begin{prop}
The map $\cP_2[d]$ is generically injective.
In other words, we can recover the element $(C_0, \beta_1, \dots, \beta_k)\in \mathcal{D}_k$ from the general element $(P, \tau)$ lying in the image of $\cP_g[d]$.

\begin{proof}
As in \cite{nos}, we consider isomorphisms $\varphi: JN\times JN\ra P$, where $N$ is a smooth genus $k(g-1)$ curve, such that \begin{equation}\tag{$\star$}
	\hat \varphi\circ \lambda_\tau\circ \varphi= \begin{pmatrix} 
		2\lambda_0 & \lambda_0\gamma\\
		\lambda_0\gamma &2\lambda_0
	\end{pmatrix},
\end{equation} for a certain $\gamma\in \text{Aut}(JN)$, and such that \begin{equation}\tag{$\star\star$}\sigma^i\circ \varphi= \varphi\circ \begin{pmatrix}
	0 &-1\\
	1&\beta_i
\end{pmatrix}, \end{equation} 
for a certain $i\in \{1,\dots, k\}$. That is, we are considering the all automorphisms $\phi_i$ of Proposition \ref{Properties}, plus possibly other ones with the same behavior. Namely, we are defining the set \[\Lambda(P, \tau,\langle \sigma \rangle):=\{(N, \varphi): \varphi: JN\times JN \cong P \text{ satisfies }(\star), (\star\star) \text{ for the same }\gamma\in \text{Aut}(JN)\}.\]

Our aim is to prove that for a generic $(P, \tau,\langle \sigma \rangle)$ in $\text{Im} \cP_{g}[d]$ we have the equality $$\Lambda(P, \tau,\langle \sigma \rangle)=\{(C_0, \beta_i), i=1, \dots, k\}.$$

Let $(N,\varphi)\in \Lambda(P, \tau,\langle \sigma \rangle) $ and let $i$ be the corresponding exponent $i$ of property $(\star\star)$. 
The composition 
\[
F: JN\times JN \xrightarrow{\varphi}P \xrightarrow{\phi_i^{-1}}JC_0\times JC_0
\] 
gives an isomorphism that pulls back the polarization of $JC_0\times JC_0$  to the one of $JN\times JN$, i.e. \[\begin{pmatrix} 
			2\lambda_0 & \lambda_0\gamma\\
			\lambda_0\gamma &2\lambda_0
		\end{pmatrix}= \hat F\circ \begin{pmatrix} 
			2\lambda_0 & \lambda_0\beta_i\\
			\lambda_0\beta_i &2\lambda_0
					\end{pmatrix}\circ F.\] 
As in \cite{nos} this gives immediately that $N\cong C_0$.


Let $ \begin{pmatrix}
	A & B\\
	C&D
\end{pmatrix}$ be the matrix associated with $F$. It is easy to see that $A^*(\lambda_{C_0})=\lambda_{C_0}$ and that an analogous argument works for $D$. Thus, $A$ and $ D$ are automorphisms of $JC_0$ compatible with the principal polarization. Hence, they come from automorphisms of $C_0$.

On the other hand, the property $(\star\star)$ determines the following equality
\[\begin{pmatrix}
	0 &-1\\
	1&\beta_i
\end{pmatrix} \circ F= F\circ \begin{pmatrix}
0 &-1\\
1&\gamma
\end{pmatrix}. \]
  An easy computation shows that it yields \[\gamma= D^{-1}\beta_iD.\]
  By Theorem \ref{teo simplicity}, the automorphism $D$ commutes with any $\beta_i$. Hence, $\gamma =\beta_i$ and we are done. 
	\end{proof}
\end{prop}

The next step is to read from $(C_0, \beta_1, \dots, \beta_k)\in \mathcal{D}_k$ the correct information to obtain the element $(C, \eta)\in \cRH_{g}[d]$. 
For this, we can apply a similar the argument to that in \cite[Section 3, Step 4]{nos}, with some minor changes. First, we observe that the map  $h_0 : C_0\ra \mathbb P^1$ determine $(C,\eta)$.
Then we apply theta-duality construction to the curves $\beta_i(C_0)$, that is, we look for the set of translates of the theta divisor that contain each of these curves. These sets determine the map $h_0$ under some numerical restrictions. To be more precise: we fix a point $x \in C_0$, and we consider the  injection $\iota_x: C_0\hookrightarrow JC_0, p\mapsto [p-x]$. The canonical theta divisor $\Theta^{can}$ in $\Pic^{k(g-1)-1}(C_0)$ is, by Riemann's parametrization Theorem, the Brill-Noether locus $W_{k(g-1)-1}(C_0)$. We define:  
\[
T'_i:=\{\xi \in \Pic^{k(g-1)-1}(C_0) \mid \beta_i(C_0)+\xi \subset \Theta^{can}= W_{k(g-1)-1}(C_0)\}.
\]
 Observe that to use the definition of $\beta_i$, we need to see $JC_0=\pi_0^*(JC_0)$ as a subvariety of $P$. Given a point $p$ in $C_{0}$, $\iota_x(p)=[p-x]$ appears as $[p'+j(p')-x'-j(x')]\in P$, where $x',p'\in \widetilde C$ are preimages of $x$ and $p$, respectively. 
     We denote by 
     \[
     p', p'_1:=\sigma (p'),\ldots ,p'_{d-1}=p'_{2k}:=\sigma ^{d-1}(p'),
     \]
the whole fiber $f^{-1}(f(p'))$, and analogously for $x'$. We denote by $p_i$ (resp. $x_i$) the image of $p'_i$ (resp. $x'_i$) in $C_0$.

\begin{prop} Assume that $(d-1)(g-1)\ge 7$ and that $g\ge 2$ is not congruent with $3$ modulo $d$. Then $T_i'=x_i+x_{d-i}+W_{k(g-1)-3}(C_0)$.
\end{prop}
Notice that the assumption $(d-1)(g-1)\ge 7$ is equivalent to $\dim W_{k(g-1)-3}(C_0)>0$.
 \begin{proof}
As in the proof of \cite[Proposition 3.13]{nos} we compute the action of $\beta_i=\sigma^i+\sigma^{-i}$ on $[p'+j(p')-x'-j(x')]$ to obtain that, as elements in $JC_0$:
     \[
        \beta_i ( [p-x])=[p_i+p_{d-i}-x_i-x_{d-i}].
     \]
     Observe that these points describe  the fiber $h_0^{-1}(h_0(x))$, more precisely, as divisors:
      \begin{equation}\label{fibre_h_0}
      h_0^{-1}(h_0(x))=x+x_1+\ldots +x_{d-1}.
      \end{equation}
By definition, $\xi \in T'_i  $ means that:
     \[
     h^0(C_0,\xi + p_i+p_{d-i}-x_i-x_{d-i})>0,\qquad \text{for all } p\in C_0.
     \]
     If $\xi $ is of the form $x_i+x_{d-i}+E$ for some effective divisor $E$ the condition is satisfied trivially, hence $x_i+x_{d-i}+W_{k(g-1)-3}(C_0)\subset  T'_i.$ 
 To prove the opposite inclusion, we consider $\xi\in T'_i$.  If $h^0(C_0, \xi-x_i-x_{d-i})>0$ then we are done. So assume $h^0(C_0, \xi-x_i-x_{d-i})=0$. Set $L:= K_{C_0}-(\xi-x_i-x_{d-i})$. Arguing as in loc. cit. we obtain that $L$ gives a $g^1_{k(g-1)+1}$  whose associated map $\varphi_L$ satisfies that $h_0^{-1}(h_0(p))$  is contained in the fiber of $\varphi_L$. Hence, $\varphi_L$ factorizes through $h_0$ and we have that $\deg(h_0)=d$ divides $\deg(L)=k(g-1)+1$. For $g=2$ this is impossible since $k+1<d=2k+1$. Assume that $g$ is at least $3$ and put $dr=k(g-1)+1=\frac {d-1}2 (g-1)+ 1$. Therefore $2rd=(d-1)(g-1)+2$ and then $g$ is congruent with $3$ modulo $d$, a contradiction. 
\end{proof}

\begin{rem}
The Proposition allows recovering intrinsically the class of the divisors $x_i+x_{d-i}$ since there are no translations leaving invariant $W_{k(g-1)-3}(C_0)$. It may well happen that $x_i+x_{d-i}$  belong to a $g^1_2$ linear series, since we have not excluded the possibility of $C_0$ being hyperelliptic. 
Assume that this is the case for a generic point in $C_0$, that is, for a generic $p$ there is an index $i$ such that $h^0(C,\cO(p_i+p_{d-i}))=2$. We can assume that $h^0(C,\cO(x_1+x_{d-i}))=2$ for the fixed point $x$ we used to embed the curve. Then $\beta_1(p-x)=p_1+p_{d-1}-x_i-x_{d-1}$. Replacing $p$ for a convenient point $p'$ in the same fiber, we obtain $\beta_1(p'-x)=p_i+p_{d-i}-x_i-x_{d-1}=0$ since both divisors represent the hyperelliptic linear series. This contradicts that $\beta_1$ is an automorphism.
\end{rem}

Since we have recovered the whole fiber $h_0^{-1}(h_0(x))$ for a generic point in $C_0$ and the curve is smooth, the whole map $h_0$ is obtained, and therefore we have finished the proof of Theorem \ref{teo Prym}. 

\begin{rem}\label{high_genus}
According to the Castelnuovo-Severi inequality \cite{ac}, the curve $C_0$ cannot admit more
than two maps to $\PP^1$ of degree $d$ if its genus satisfies
$g(C_0)= (g-1)(d-1)> (d-1)^2$, that is if $g>2d-1$. This implies that even in this range the map $\cP_2[d]$ is injective. Therefore, one can conclude that the generic injectivity of $\cP_g[d]$ holds as well.
\end{rem}

\end{document}